\documentclass[11pt]{article}
\usepackage{amsfonts}      %para las doublestroke
\usepackage[pdftex]{graphicx}
\usepackage[margin=3cm]{geometry}
\usepackage{graphicx}
\usepackage{graphics}
\usepackage{amssymb,amsmath}
\usepackage{enumerate}
\usepackage{amscd}
\usepackage{bibentry}
\usepackage[sort]{natbib}
\nobibliography*
\usepackage{subfigure}
\usepackage{fancyhdr}
\usepackage{amsthm}
\usepackage{color}

\usepackage[utf8]{inputenc}

\usepackage{fancyheadings}

\usepackage{bigfoot}
\usepackage{alphalph}

\DeclareNewFootnote{M}[alph] % My footnotes
\MakePerPage{footnoteM}
\DeclareNewFootnote{F}         % De Finetti's Original Footnotes
\MakePerPage{footnoteF}

\newtheorem{lemma*}{Lemma}

\makeatletter
\renewcommand*\env@matrix[1][*\c@MaxMatrixCols c]{%
  \hskip -\arraycolsep
  \let\@ifnextchar\new@ifnextchar
  \array{#1}}
\makeatother

\usepackage{cleveref}
\crefname{section}{§}{§§}
\Crefname{section}{§}{§§}

\title{A translation of ``The characteristic function of a random phenomenon'' by Bruno de Finetti.\footnote{Citations of this translation should refer also to de Finetti’s original paper, as cited in the Abstract.}}
\date{\today}
\usepackage{authblk}
\author[$\ddagger$]{David Alvarez-Melis\footnote{Corresponding author. Email: \texttt{davidam@csail.mit.edu}.}}
\author[$\ddagger$]{Tamara Broderick}
\affil[$\ddagger$]{ Computer Science and Artificial Intelligence Lab, MIT }

\begin{document}

\maketitle

% \pagestyle{fancy}
% \lhead{B. de Finetti}
% \rhead{Translation by: Alvarez-Melis, Broderick}

\vspace{-0.7cm}
\abstract{This article is a translation of Bruno de Finetti's paper ``Funzione Caratteristica di un fenomeno aleatorio'' which appeared in \textit{Atti del Congresso Internazionale dei Matematici, Bologna 3-10 Settembre 1928, Tomo VI, pp.~179-190}, originally published by Nicola Zanichelli Editore S.p.A. The translation was made as close as possible to the original in form and style, except for apparent mistakes found in the original document, which were corrected and are mentioned as footnotes. Most of these were resolved by comparing against a longer version of this work by de Finetti, published shortly after this one under the same title\footnoteM{Published in 1930 as part of the Memories of the Royal Lincean Academy. Full details are provided by de Finetti himself is his second footnote in this page.}. The interested reader is highly encouraged to consult this other version for a more detailed treatment of the topics covered here. Footnotes regarding the translation are labeled with letters to distinguish them from de Finetti's original footnotes.
% Comunicazioni Sezione IV [A] - V - VII

\section*{The original document}

% =================  Section  1 ==================%
\indent  § 1. - The purpose of this report is to show how the characteristic function method\footnoteF{V. i Trattati di \textit{Calcolo delle Probabilita} di G.\ Castelnuovo e di P.\ Levy.}, already successfully introduced in the theory of random variables, is also usefully suitable for the study of random phenomena. We will therefore show, in the first place, how a random phenomenon can be completely determined by means of its characteristic function, and we will then sketch the operations that make the latter a powerful computational tool. A longer account will be found in an essay\footnoteM{Lit.~\textit{memoria}. The author is referring to the longer version of the paper, as explained in the Abstract.} that will be presented at the earliest to the Royal Lincean Academy\footnoteF{Memorie della R. Acc. Naz. dei Lincei, S.$6^a$, vol. IV, fasc. V.}.

% \footnote{Lit. \textit{memoria}. The essay the author is referring to is a longer, more detailed version of this work, published one year later and with the same title. The interested reader is highly encouraged to consult this other version for the sake of completion.}

%  which the reader is advised to consult: B. de Finetti, Funzione Caratteristica di un Fenomeno Aleatorio, in \textit{Memmorie della R. Accademia dei Lincei}, 1930.]

\vspace{0.5cm}
§ 2. - We will refer to a phenomenon of which any number of trials can be carried out (or can at least be conceived) as a \textit{random phenomenon} when the order in which favorable and unfavorable outcomes alternate can only be attributed to chance. That is, we require that all $\binom{n}{h}$ sequences of $n$ trials with $h$ successes - which differ only in the order of the events - have the same probability. This, in specific terms, is the characteristic property of the phenomena which we have convened to define as random.

It will be useful to show this restriction with an example, and thus have a clear idea of the topic of this work. If we have a coin or a die and we always throw it in the same way, there is no reason, of a causal type, not even if we are not sure about the faultlessness of the coin, which can influence the order in which favorable and unfavorable trials alternate: the order will be determined by chance, and thus we have a random phenomenon according to the definition above. The same can be said for the game of roulette, for the drawing from an urn whose elements have been chosen from a known collection, and all similar cases. If instead we consider a sequence of shots on a target by the same shooter, or the succession of rainy and non rainy days, or the days in which the neighbor\footnoteM{Lit.~\textit{signore di rimpetto}.} shaves his beard, such condition cannot be reasonably believed to hold. This is because in the first example we can on the one hand foresee a progressive training of the shooter, but on the other hand also an increasing fatigue, which makes likely an accumulation of favorable results in a short period of time. In the second example, days with rain will be accumulated in longer or shorter rainy periods, not to mention seasonal periodicity, and in the last example the neighbor will always shave at relatively regular intervals. 

To decide, in practice, if a given phenomenon can be considered to be random or not it suffices to think whether an eventual regularity or another singularity found in the order of the sequence could be attributed to chance (therefore having a random phenomenon) or whether this could be attributed to some circumstance related to the phenomenon, so that it would be conceivable that in another equivalent sequence of trials the same order would likely be repeated. 

\vspace{0.5cm}
§ 3. - If $n$ trials of a given random phenomenon are carried out, the number of those resulting in a favorable outcome is clearly a random variable $x_n$ that can take only the values $0,1,\dots,n$. If we denote by $\omega_h^{(n)}$ the probability that the phenomenon in question is verified\footnoteM{The text reads \textit{verificare}. The closest translation was kept, although a more natural way to say this in English would be ``results in a success''.} $h$ times out of $n$ trials, the random variable $x_n$ is characterized by the probabilities $\omega_0^{(n)}, \omega_1^{(n)}, \dots, \omega_n^{(n)}$ of occurrence of all possible outcomes. In the particular and well known case where the phenomenon has a constant probability $p$ known in advance, we know that $\omega_h^{(n)} = \binom{n}{h} p^h(1-p)^{n-h}$, but in the general case which concerns us the $\omega_h^{(n)}$ can be of an arbitrary form (except for the restrictions posed by the nature of the problem, which we address below).

A random phenomenon will thus define a sequence of random variables $x_1, x_2, \dots, x_n, \dots,$ which must naturally turn out to be dependent of each other. Such dependence translates analytically into a recurrent differential relation that links their characteristic functions $\psi_1,\psi_2\dots \psi_n,\dots$, and which constitutes the foundation of this work. 
We prove that as $n$ grows the function
\[ \psi_n\left ( \frac{t}{n} \right) = {\sum_{h=0}^n} \omega_h^{(n)}e^{i \frac{h}{n}t} \]
tends uniformly in every finite region to the entire function
\[ \psi(t) = {\sum_{h=0}^{\infty}} \omega_h^{(h)} \frac{i^ht^h}{h!}, \]
which is precisely what we will define as the \textit{characteristic function of the random phenomenon}. Given $\psi$, all the $\psi_n$ - and consequently, all the $\omega_h^{(n)}$ - can be obtained, which justifies the name of this function. 

The integral from $-\infty$ to $\infty$ of the function $\frac{e^{it}-e^{-i \xi t}}{it}\psi(t)$ exists for every value of $\xi$, and equals $0$ for $\xi <0$ and $2\pi$ for $\xi >1$. Consequently, there exists a random variable of which $\psi(t)$ is the characteristic function, with corresponding cumulative distribution function given by
\[ \Phi(\xi) = \frac{1}{2\pi} \int_{-\infty}^{\infty} \frac{e^{it}-e^{-i \xi t}}{it}\psi(t)dt, \] 
with $\Phi(\xi) =0$ for $\xi <0$ and $\Phi(\xi) =1$ for $\xi >1$.

From these results two important theorems follow: \\
\indent I. The probability that the frequency over $n$ trials is contained within given limits $\xi_1$ and $\xi_2$ tends to $\Phi(\xi_2) - \Phi(\xi_1)$ as $n$ grows. \\
\indent II. The probability that all the frequencies after the $n$-th one are contained within given limits $\xi_1$ and $\xi_2$ tends to\footnoteM{The author uses the notation limd (\textit{limite destro} - right limit) and lims (\textit{limite sinistro} - left limit) for these one-sided limits, but subtracts the latter from the former, which considering the interval of interest and comparing to the longer version of the paper, is incorrect. We show here the correct order of the limits and use the conventional one-sided notation.}
\[    \lim_{\xi \rightarrow \xi_2^-}  \Phi(\xi) - \lim_{\xi \rightarrow \xi_1^+} \Phi(\xi) \]
as $n$ grows.

Given a $\psi(t)$, in order for a random phenomenon having this characteristic function to exist, it is necessary and sufficient that the CDF (naturally real and never decreasing) vanish for $\xi<0$ and equal $1$ for $\xi >1$. 

\vspace{0.5cm}

§ 4. - We briefly present the calculations.

 The following relation must hold between the $\omega_h^{(n)}$:
 \begin{equation}
 	\omega_k^{(m)} = \sum_{h = k}^{n-m+k} \omega_{h}^{(n)} \frac{\binom{h}{k}\binom{n-h}{m-k}}{\binom{n}{m}}
 \end{equation}
because $\binom{h}{k}\binom{n-h}{m-k}/\binom{n}{m}$ is the probability that $k$ out of $m$ trials are successful, taken from $n$ of which $h$ are successful, when all the combinations are equally likely. In particular (for $m=n-1$):
\begin{equation}\label{due}
	n \omega_{k}^{(n-1)} = (n-k)\omega_{k}^{(n)} + (k+1)\omega_{k+1}^{(n)}
\end{equation}
and setting
\begin{equation}\label{tre}
	\Omega_n(z) = \sum_{h=0}^n \omega_h^{(n)}z^h
\end{equation}
all the equations in \eqref{due} for $k=0,1,\dots,n-1$ are summarized in the differential recurrence relation
\begin{equation}\label{quattro}
	n\Omega_{n-1}(z) = n\Omega_n(z) + (1-z)D\Omega_n(z).
\end{equation}
Taking derivatives in \eqref{quattro} and obtaining the value of the consecutive derivative for $z=1$ we obtain
\begin{equation}\label{cinque}
	\Omega_n(1+z) = \sum_{h=0}^n \binom{n}{h}\omega_{h}^{(h)}z^h
\end{equation}
and we find that as $n \rightarrow \infty$, $\Omega_n(1 + \frac{z}{n})$ tends uniformly to 
\begin{equation}\label{sei}
	\Omega(1+z) = \sum_{h=0}^{\infty} \omega_{h}^{(h)}\frac{z^h}{h!}.
\end{equation}
 The characteristic function $\psi_n$ is 
\begin{equation}\label{sette}
	\psi_n(t) = \Omega_n(e^{it}).
\end{equation}
It can be proven that\footnoteM{The original reads $\Omega(e^{\frac{t}{n}})$, which is incorrect. See the equivalent equation in the long version of this essay.} $\Omega_n \left( e^{\frac{t}{n}} \right)$ in turn tends to $\Omega(1+t)$, and therefore the characteristic function of the random phenomenon (the function to which $\psi_n(\frac{t}{n})$ tends uniformly) is:
\begin{equation}\label{otto}
	\psi(t) = \Omega(1 + it).
\end{equation}

\vspace{0.5cm}

§ 5. - A first noteworthy consequence: \eqref{otto} and \eqref{sei} show that as $n$ goes to infinity, the $m$-th moment of $\frac{x_n}{n}$, namely the limit-probable-value\footnoteM{Lit.~\textit{valor-probabile-limite}, presumably referring to the expected value.} of the $m$-th power of the frequency for arbitrarily large number of trials, tends to $\omega_m^{(m)}$, that is, the probability that all $m$ trials are favorable. 

The theorems stated in § 3 are obtained by translating the results obtained for the characteristic function into those corresponding to the CDF. Denoting by $\Phi_n(x)$ the CDF of $x_n$, we have:
\begin{equation}\label{nove}
	\lim_{n \rightarrow \infty} \Phi_n(n\xi) = \Phi(\xi).
\end{equation}
The other theorem\footnoteM{The author is referring to Theorem II of Section 3, as can be confirmed in the long version of the paper.}, relating to the probability that all the frequencies beyond a certain $n$ belong to a given interval, can be thought of as a generalization to the case of an arbitrary random phenomenon of an analogous result for independent variables with equal probability. The latter would be a particular case of the well-known {\sc Cantelli} Theorem\footnoteF{F.P.\ Cantelli: \textit{Sulla probabilit\`a come limite della frequenza.} Rend.\ R.\ Acc.\ dei Lincei, serie V, vol.\ XXVI, gennaio 1917.}\footnoteM{The author presumably refers to what is known today as the Second Borel-Cantelli Lemma.}. The proof of the theorem stated here can also be brought back to the case addressed by {\sc Cantelli}.

\vspace{0.5cm}

§ 6. -  Equations \eqref{sei}, \eqref{cinque} and \eqref{sette} clearly prove the claim that $\psi$ suffices to completely determine all the $\psi_n$. The same can be said for $\Phi$, because given this function we have
\begin{gather}
	\psi(t) = \int_0^1 e^{it\xi}d\Phi = e^{it} - it\int_0^1 e^{it\xi}\Phi(\xi)d\xi \\ \label{dieci}
	\Omega_n(1 + z) = \int_0^1 (1 + z\xi)^n d\Phi = (1+z)^n -nz\int_0^1(1+z\xi)^{n-1}\Phi(\xi)d\xi \\
	\omega_h^{(n)} = \binom{n}{h} \int_o^1 \xi^h(1-\xi)^{n-h} d\Phi = \binom{n}{h}\int_0^1(h -n\xi)\xi^{n-1}(1-\xi)^{n-h-1}\Phi(\xi)d\xi.
\end{gather}

\vspace{0.5cm}

§ 7. -  We consider two cases of particular importance, which will ultimately serve as examples too.  

In the well-known case where the probability of a phenomenon is known a priori to be equal to $p$,
\begin{gather*}
\omega_h^{(n)} = \binom{n}{h}p^h(1-p)^{n-h}, \qquad \omega_h^{(h)} = p^h, \\
\Omega_n(1+z) = (1+pz)^n, \qquad \Omega(1+z) = \lim_{n\rightarrow \infty} \left(1 + \frac{pz}{n}\right )^n = e^{pz},\\
\psi(t) = e^{ipt},
\end{gather*}
or, equivalently (directly from \eqref{sette}):
\[ \psi(t) = \sum_{h=0}^{\infty} p^h \frac{i^ht^h}{h!} = e^{ipt}, \qquad \Phi(\xi)= 
\begin{cases}
	0 \mbox{ if }& \xi<p \\
	\frac{1}{2} \mbox{ if }& \xi =p \\
	1 \mbox{ if }& \xi >p \\
\end{cases} \]
As $n$ grows the probability that the frequency is contained in a neighborhood $p\pm \epsilon$ of $p$ tends to one; from this it follows that\footnoteM{The original reads \textit{inversamente da tale ipotesi scende}, literally translated as ``inversely from this hypothesis it follows.'' Based on the preceding and following sentences, it seems that the author is simply using this transition as a way to introduce a cause-effect statement, hence the translation used here.} $\psi(t) = e^{ipt}$, and consequently, that the phenomenon has probability constant and equal to $p$ in all the trials. As a particular case, for $p=0,p=1$ we have $\psi(t) =1, \psi(t)=e^{it}$.
In the case where all the possible values for the frequencies are equally likely we will have
\begin{gather*}
	\omega_0^{(n)} = \omega_1^{(n)} = \dots = \omega_{n}^{(n)} = \frac{1}{n+1}, \quad \omega_h^{(h)} = \frac{1}{h+1}, \\
	\Omega_n(z) = \frac{1}{n+1} (1 + z + z^2 + \dots + z^n) = \frac{1}{n+1} \cdot \frac{1 - z^{n+1}}{1-z}, \\
	\psi_n(t)= \frac{1}{n+1}\cdot \frac{1 - e^{it(n+1)}}{1-e^{it}}, \\
	\Omega(1+z) = \lim_{n\rightarrow \infty} \frac{1}{n+1} \cdot \frac{\left( 1+ \frac{z}{n}\right)^{n+1} -1}{ \frac{z}{n}} = \frac{e^z -1 }{z}, \\
	\psi(t) = \frac{e^{it}-1}{it} = \sum_{h=0}^{\infty} \frac{t^h}{(h+1)!}, \\
	\Phi_n(h) = \frac{h + \frac{1}{2}}{n+1} \quad (h=0,1,\dots,n), \qquad \Phi(\xi) = \lim_{n \rightarrow \infty} \Phi_n(n\xi) = \xi.
\end{gather*}
As $n$ grows, the probability that the frequency is contained in an interval between $\xi_1$ and $\xi_2$ tends to $\xi_2 -\xi_1$; from here it follows that\footnoteM{See previous note.}
\[ \psi(t) = \int_0^1 e^{i\xi t}d\xi = \left [ \frac{e^{it\xi}}{it}  \right]_0^1 = \frac{e^{it} -1}{it}, \qquad \omega_h^{(n)} = \frac{1}{n+1}, \]
and thus the phenomenon has equal probability among all possible frequencies over $n$ trials.

\vspace{0.5cm}

§ 8. - We now move on to the operations on the characteristic functions. 

As a general observation we can say that all the operations that we will encounter are distributive, up to a multiplicative factor (whenever necessary) which makes $\psi(t)$ have the value $1$ for $t=0$, as must necessarily happen. 

By introducing the operator $U$:
\[ Uf(t) = \frac{f(t)}{f(0)} \]
we can say that the operations that will be presented are products of the form $UF$ with $F$ distributive.

Given that the CDF $\Phi(\xi)$ is a bijective linear function of the characteristic function $\psi(t)$, to any distributive operation on $\psi$ will correspond the transformation that operates on $\Phi$.

Two operations useful for simplifying the notation are $\mathbf{P}_n$ (read: \textit{the n-th polynomial}) that applied to $\psi$ yields $\Omega_n$\footnoteF{Precisely the function $\Omega_n(1+z)$, not $\Omega_n(z)$.}, and ${ h \brack n}$, which applied to $\psi$ yields $\omega_h^{(n)}$. They can be defined in general as follows. If
\[ f(t) = \sum_{h=0}^{\infty} a_h \frac{i^ht^h}{h!},  \]
let
\begin{equation}\label{tredici}
	\mathbf{P}_nf(t) = \sum_{h=0}^{n} \binom{n}{h} a_ht^h = \sum_{h=0}^{n} (1+t)^h \left \{ { h \brack n} f \right \}.
\end{equation}
Note in particular that $ { n \brack n} f = a_n$, and that $\psi_n(t) = \mathbf{P}_n \psi(e^{it} -1 )$.

\vspace{0.5cm}

§ 9. -  Let $\psi(t)$ be the characteristic function of some random phenomenon. The characteristic function of the complementary event is
\begin{equation}\label{quattordici}
	K \psi(t) = e^{it}\psi(-t).
\end{equation}
In fact, saying that $h$ out of $n$ trials are successes is equivalent to saying that $n-h$ are successes for the complementary event, which is expressed by
\begin{equation}\label{quindici}
	{h \brack n} K = { n-h \brack n},
\end{equation}
and yields (writing $\omega_h^{(n)} = {n \brack h } \psi$):
\begin{gather*}
	\begin{split}
   	 \mathbf{P}_nK\psi(z-1) = \omega_n^{(n)} + \omega_{n-1}^{(n)}z + \dots + \omega_0^{(n)}z^n = \phantom{a+a+a+a+a+a+a+a} \\
   	  = z^n \{ \omega_0^{(n)} + \omega_1^{(n)}\frac{1}{z} + \dots + \omega_n^{(n)} \frac{1}{z^n} \} = z^n\mathbf{P}_n \psi\left(\frac{1}{z} -1 \right),
	\end{split}
	\\
	  \mathbf{P}_n K \psi(e^{it} -1) = e^{nit} \mathbf{P}_n \psi(e^{-it} -1)
\end{gather*}
and
\[ K \psi(t) = \lim_{n \rightarrow \infty} \mathbf{P}_n K \psi( e^{i \frac{t}{n}} -1) = \lim_{n\rightarrow \infty} ( e^{i \frac{t}{n}} )^n \lim_{n \rightarrow \infty} \mathbf{P}_n \psi( e^{-i \frac{t}{n}} - 1) = e^{it} \psi(-t). \]
In particular $ {n \brack n } K = { 0\brack n}$, so
\begin{equation}\label{sedici}
	K \psi(t) = 1 + \omega_0^{(1)}it - \omega_0^{(2)} \frac{t^2}{2!} - \omega_0^{(3)}i \frac{t^3}{3!} + \dots + \omega_0^{(n)} i^n \frac{t^n}{n!} + \dots
\end{equation}
If $\Phi$ is the CDF corresponding to $\psi$, then we find that the CDF $K_{\Phi}\Phi$ corresponding to $K\psi$ is
\begin{equation}\label{diciassette}
	 K_{\Phi} \Phi(\xi) = 1 - \Phi(1 - \xi)
\end{equation}
which also follows intuitively: the probability that in the limit the frequency of an event is less than $\xi$ is equal to the probability that the frequency of the complementary event is greater than $1- \xi$.

$K$ is a distributive, reversible and involutory function:
\[ K(\psi' + \psi'') = K\psi' + K \psi'', \qquad KK\psi = \psi, \qquad K^{-1}\psi= K\psi. \]

\vspace{0.5cm}

§ 10. - Under the hypothesis that the first trial is a success, the characteristic function is 
\begin{equation}\label{diciotto}
	R \psi = UD \psi,
\end{equation}
while under the hypothesis that the first trial is a failure, it is
\begin{equation}\label{dicianove}
	S \psi  = U(i - D)\psi
\end{equation}
In general, after $r$ successes and $s$ failures, the characteristic function becomes\footnoteM{The original text has $(1-D)$ instead of $(i-D)$ in \eqref{venti}, which is incorrect. Compare with the equivalent equation in the longer version of the essay.}
\begin{equation}\label{venti}
	R^rS^s \psi = UD^r(i-D)^s \psi.
\end{equation}

The probability $\omega_n^{(n)}$ that the first $n$ trials are all favorable is in fact equal to the product of $\omega_1^{(1)}$, the probability that the first trial is successful, times the probability that, after this hypothesis being confirmed, the following $n-1$ trials are also favorable. This probability is ${n-1 \brack n-1}R\psi$. The $n$-th coefficient in the expansion of $R\psi$ is therefore $\omega_{n+1}^{(n+1)}$, the $(n+1)$-th coefficient in the expansion of $\psi$, divided by the first coefficient:
\[ \omega_1^{(1)} = -i D \psi(0), \]
and thus
\[ R\psi(t) = \frac{1}{\omega_1^{(1)}} \left \{ \omega_1^{(1)} + \omega_2^{(2)}\frac{it}{2!} - \omega_3^{(3)}\frac{t^2}{2!} + \dots + \omega_{n+1}^{(n+1)}\frac{i^nt^n}{n!}  + \dots  \right \} = \frac{-iD\psi(t)}{-iD\psi(0)} = UD\psi(t). \]
To show that $S = U(i-D)$ it is convenient to start from the observation that $S$ is clearly a transformation of $R$ under $K$: $S = KRK$. From \eqref{quattordici}:
\begin{gather*}
	DK\psi(t) = e^{it}[i\psi(-t) - D \psi(-t)]; \\
	KDK \psi(t) = e^{it} \{ e^{-it}[i\psi(t) - D \psi(t)]\} = (i-D)\psi(t); \qquad KDK = i -D;
\end{gather*}
and by the properties of $U$:
\[ S  = KRK = UKDK = U(i-D), \]
from which 
\[ S \psi(t) = \frac{i \psi(t) - D \psi(t)}{i  - D \psi(0)} = \frac{1}{1 - \omega_1^{(1)}} \sum_{h=0}^{\infty} ( \omega_h^{(h)} - \omega_{h+1}^{(h+1)}) \frac{i^ht^h}{h!}. \]

The operators $R$ and $S$ commute:
\[ RS = SR,\] 
hence
\[ R^r S^s = S^s R^r = S^{s_1} R^{r_1}S^{s_2} R^{r_2} \dots \quad (s_1 + s_2 + \dots = s, \medspace r_1 + r_2 + \dots = r).\]
This proves that after $r$ successes and $s$ failures, independently of the order in which they occur, the characteristic function becomes $R^rS^s \psi$.

For a random phenomenon with characteristic function $\psi$, and under the hypothesis that the first $r+s$ trials\footnoteM{The original text reads $r+1$, but given the following sentence, the total number of trials must be $r+s$, as written here.} result in $r$ favorable and $s$ unfavorable outcomes, the probability that a total of $h$ out of $n$ trials will be successful is given by the formula
\begin{equation}\label{ventuno}
		{ h \brack n} R^rS^s\psi = \frac{\binom{h+r}{r} \binom{n-h+s}{s}}{\binom{n+r+s}{n}} \cdot \frac{ { h + r \brack n+r+s } \psi }{ { r  \brack r+s} \psi  }\medspace.
\end{equation}

\vspace{0.7cm}

§ 11. -  We denote by $R_{\Phi}, S_{\Phi}$ the operations that transform the CDF corresponding to $\psi$ to those corresponding to $R \psi$ or to $S \psi$ respectively.
Starting from the expression of $D\psi(t)$ obtained from differentiating \eqref{dieci} we obtain
\begin{equation}\label{ventidue}
	R_{\Phi}\Phi(\xi) = \frac{ \xi \Phi(\xi) - \int_0^{\xi}\Phi(\lambda)d\lambda }{1 - \int_0^1\Phi(\lambda)d\lambda} = \frac{ \int_0^{\Phi(\xi)}\xi d\Phi }{\int_0^1 \xi d\Phi};
\end{equation}
analogously, we see that 
\begin{equation}
	S_{\Phi}\Phi(\xi) = \frac{\int_0^{\Phi(\xi)}(1-\xi)d\Phi}{\int_0^1(1-\xi)d\Phi},
\end{equation}
and in general
\begin{equation}\label{ventiquattro}
	R_{\Phi}^rS_{\Phi}^s\Phi(\xi) = \frac{\int_0^{\Phi(\xi)}\xi^r(1-\xi)^sd\Phi}{\int_0^1\xi^r(1-\xi)^sd\Phi}.
\end{equation}

We can give a more expressive form to these results. By the mean value theorem:
\begin{equation}\label{venticinque}
	[R_{\Phi}^rS_{\Phi}^s\Phi]_{\xi_1}^{\xi_2} = \frac{\bar{\xi}^r(1-\bar{\xi})^s}{\int_0^1\xi^r(1-\xi)^sd\Phi}[\Phi]_{\xi_1}^{\xi_2}
\end{equation}
with $\xi_1 \leq \bar{\xi} \leq \xi_2$ and denoting $[\Phi]_{\xi_1}^{\xi_1}:= \Phi(\xi_2) - \Phi(\xi_1)$. 

From this formula follows a remarkable asymptotic theorem from the domain of relations between probability and frequency. If the frequency over a sufficiently large number of trials is $f$, the characteristic function tends to the one corresponding to the case where the phenomenon had known prior probability $f$, unless $\Phi$ is constant in a neighborhood of $f$. More precisely, if for no $\epsilon >0$ the following holds 
\[ \Phi \left(\frac{r}{r+s} - \epsilon \right) = \Phi \left(\frac{r}{r+s} + \epsilon \right) \]
then 
\[ \lim_{n\rightarrow \infty} (R_{\Phi}^rS_{\Phi}^s)^n \Phi(\xi) = \begin{cases}
	0 \mbox{ if }& \xi < f \\
	\frac{1}{2} \mbox{ if }& \xi = f \\
	1 \mbox{ if }& \xi > f \\
\end{cases} \]
whence it follows that
\begin{equation}\label{ventisei}
	\lim_{n \rightarrow \infty} (R^rS^s)^n \psi(t) = e^{ift},
\end{equation}
where the convergence is uniform in any finite region ({\sc Castelnuovo}, Op.\ cit.,\ vol.\ II, p.\ 198).

Therefore, whatever the nature of a random phenomenon is (as long as its function $\Phi$ satisfies the restrictions mentioned previously), the probability that after $n$ trials with frequency $f$ the frequency on the following trials tends to a limit that differs from $f$ by more than a given $\epsilon$ can be made arbitrarily small by taking $n$ sufficiently large.

\vspace{0.5cm}

§ 12. -  As an exercise, we apply the results found above to the characteristic functions considered in § 8.

If $\psi(t) = e^{ipt}$, then
\begin{equation*}
	K \psi(t) = e^{i(1-p)t}, \qquad R \psi(t) = S \psi(t) = R^rS^s\psi(t) = e^{ipt}.
\end{equation*}
There are no other characteristic functions that remain invariant upon knowledge of the success of a trial: if $R\psi = \psi$ or $S\psi = \psi$, it follows that $\psi(t) = e^{ipt} \medspace (0 \leq p \leq 1)$. 

If $\psi(t) = \frac{e^{it}-1}{it}$:
\begin{gather*}
	K\psi(t) = \psi(t) = \frac{e^{it}-1}{it}; \\
	R\psi(t) = \frac{2}{t^2}(e^{it} - ite^{it} -1); \quad S\psi(t) = \frac{2}{t^2}(1 + it - e^{it});
\end{gather*}
\begin{equation}
	{ h \brack n} R^rS^s\psi = \frac{\binom{h+r}{r} \binom{n-h+s}{s}}{\binom{n+r+s}{n}} \cdot \frac{n+r+s+1}{r+s+1}.
\end{equation}
In particular, the probability that the $(r+s+1)$-nth trial is favorable after $r$ successes and $s$ failures is given by
\begin{equation}\label{ventotto}
	{ 1 \brack 1}  R^rS^s\psi = \frac{r+1}{r+s+2}.
\end{equation}
It is this formula that is frequently used (and even abused) in the theory of posterior probabilities. It is rigorously exact when $\psi(t) = \frac{e^{it} - 1}{it}$, but is valid only in this very special case. 

\vspace{0.5cm}

§ 13. -  Two other problems deserve a brief description. 

Given two phenomena (analogously for three or more), independent of each other, with characteristic functions
\[ \psi'(t) = \sum_{h=0}^{\infty} a_h\frac{i^ht^h}{h!}, \qquad \psi''(t) = \sum_{h=0}^{\infty} b_h\frac{i^ht^h}{h!},\]
the event of both having a favorable outcome is a random phenomenon with characteristic function
\begin{equation}\label{ventinove}
	\psi(t) = \sum_{h=0}^{\infty} a_hb_h\frac{i^ht^h}{h!}.
\end{equation}
Indeed, $a_h$ is the probability that the outcomes of the first phenomenon in $h$ trials are all favorable, and analogously for $b_h$, so that $a_hb_h$ is the probability that both of these phenomena result in all favorable outcomes after $h$ trials.

In particular if $\psi''(t) = e^{ipt}$ (phenomenon with known probability $p$) we have $\psi(t) = \psi'(pt)$. If 
\[ \psi''(t) = \frac{e^{it} -1}{it} \]
then 
\[ \psi(t) = \frac{1}{t} \int_0^t \psi'(t)dt.  \]
As an example, for 
\[ \psi'(t) = e^{ipt}, \qquad \psi''(t) = \frac{e^{it} -1}{it}: \qquad \psi(t) = \frac{e^{ipt} -1}{ipt}. \]

\vspace{0.5cm}

§ 14. -  If a phenomenon can depend on various mutually exclusive\footnoteM{Lit.~\textit{incompatibili}.} causes that have probabilities $\lambda_1, \lambda_2, \dots, \lambda_m$, respectively, and under these different hypotheses the phenomenon has respective characteristic functions $\psi^{(1)}(t), \psi^{(2)}(t), \dots, \psi^{(m)}(t)$, then the characteristic function of the random phenomenon is:
\begin{equation}\label{trenta}
	\psi(t) = \lambda_1\psi^{(1)}(t) + \lambda_2\psi^{(2)}(t) + \dots + \lambda_m \psi^{(m)}(t).
\end{equation}
The theory of the probability of hypotheses can be based upon this theorem in a flawlessly formal way. The classic example to which this result is applied is that of drawing from an urn that has been chosen from a given collection of urns. If we know that the proportion of black balls can be $p_1, p_2 \dots,p_m$ with probabilities $\lambda_1,\lambda_2,\dots,\lambda_m$, the characteristic function will be\footnoteM{In the original text, the last term in this equation is incorrectly given as $\lambda_m e^{i\lambda_mt}$.}:
\[ \psi(t) = \lambda_1 e^{ip_1t} + \lambda_2 e^{ip_2t} + \dots + \lambda_m e^{ip_mt}.\]
Another example - still referring to the extraction from an urn, for the sake of consolidating this idea - that better relates to the type of problems that might arise in practice is the following. Suppose we have an urn $A$ containing $n$ black and white balls, which were chosen by someone who had $N=cn$ balls ($c$ an integer greater than $1$) at his disposal, of which $H=ch$ were white and $K=ck$ were black ($H+K = N$, that is, $h + k = n$). Of all the possible hypotheses, we believe only the two following to be feasible: $a)$ the $n$ balls were chosen at random from the $N$ available, $b)$ the person who prepared the urn chose the $n$ balls in such way that the proportion of white and black balls was preserved (thus picking $h$ white and $k$ black balls). Furthermore, we know the probabilities, say $\alpha$ and $\beta$, of these two hypotheses. The event of drawing a white ball from the urn $A$ has characteristic function:
\[ \psi(t) = \alpha \cdot \psi^{(\alpha)}(t) + \beta \cdot \psi^{(\beta)} (t) \]
where 
\[ \psi^{(\beta)}(t) = e^{i \frac{h}{n}t }, \]
and
\[ \psi^{(\alpha)}(t) = \frac{1}{\binom{N}{n}} \sum_{l=0}^n \binom{H}{l} \binom{N-H}{n-l} e^{i \frac{l}{n}t} \]
(where $\binom{H}{l} \binom{N-H}{n-l}/ \binom{N}{n}$ is the probability that $l$ out of $n$ balls deposited in the urn were white). After $r+s$ draws of which $r$ yielded white balls, the characteristic function is
\[ R^rS^s\psi(t) = \frac{1}{{r \brack r+s} \psi} \left\{ \alpha \cdot { r \brack r+s} \psi^{(\alpha)}\cdot R^rS^s\psi^{(\alpha)}(t) + \beta \cdot { r \brack r+s} \psi^{(\beta)} \cdot R^rS^s\psi^{(\beta)}(t)   \right\}. \]
Taking derivatives and setting $t=0$, we obtain the probability of getting a white ball in the $(r+s+1)$-th draw (\textit{determining a posterior probability}). The probability that the balls were chosen randomly (hypothesis $a$) after having drawn $r$ white and $s$ black balls is
\[ \alpha \cdot \frac{ { r \brack r+s} \psi^{(\alpha)}}{ { r \brack r+s} \psi} \]
(\textit{determining the probability of a hypothesis}\footnoteM{The two problems described here, that of determining the posterior probability and the probability of a hypothesis, are then referred to by the author in the next section using the terms \textit{problema di probabilità a posteriori} (problem of posterior probability), and \textit{problema delle probabilità delle ipotesi} (problem of hypotheses probability), which we translate as ``hypothesis likelihood problem''.}).

For the sake of a numerical example, suppose the urn has 6 balls chosen from a total of 12 available balls, of which 4 were white and 8 were black, and it is believed that $\alpha =\frac{2}{3}$ and $\beta=\frac{1}{3}$, then\footnoteM{The denominator in the exponent of the fourth term in $\psi^{(\alpha)}(t)$ is missing in the original.}
\begin{gather*}
	\psi^{(\alpha)}(t) = \frac{1}{33} \left \{   1 + 8e^{\frac{it}{6}} + 15e^{\frac{it}{3}} + 8e^{\frac{it}{2}} + e^{\frac{2it}{3}} \right \}, \qquad \psi^{(\beta)}(t) = e^{\frac{it}{3}}, \\
	 \psi(t) = \frac{1}{99} \left \{   2 + 16e^{\frac{it}{6}} + 63e^{\frac{it}{3}} + 16e^{\frac{it}{2}} + 2e^{\frac{2it}{3}} \right \}.
\end{gather*}
After $r+s$ draws, of which $r$ yielded a white ball, the probability of hypothesis $b)$ being true (i.e.\ non random choice) is 
\[ \frac{33\cdot2^r\cdot 4^s}{16\cdot5^s + 63\cdot2^r\cdot4^s + 16\cdot 3^{r+s} + 2\cdot4^r\cdot2^s}. \] 
After 6 draws hypothesis $b)$ has probability $0.088353$ if $r=6$ white balls have been drawn, $0.176707$ if $r=5$, $0.279365$ if $r=4$, $0.359918$ if $r=3$, $0.389812$ if $r=2$, $0.353947$ if $r=1$, and finally $0.260018$ if $r=0$, that is, if only black balls have been drawn so far. As the number of trials grows, the probability of hypothesis $b)$ tends to $\frac{33}{63} = 0.523810$, $\frac{33}{79} = 0.517722$ or to zero, respectively, depending on whether the frequency is within 
\[ (\log \tfrac{5}{4})/(\log \tfrac{5}{2}) = 0.243529 \quad \text{  and  } \quad (\log \tfrac{4}{3})/(\log 2) = 0.415037,  \]
it is equal to either of these limits, or it is outside this interval.

\vspace{0.5cm}

§ 15. -  The problem of posterior probabilities consists of trying to determine the probability of a random phenomenon based on the observed frequency of successes in a given number of trials. Based on what has been noted above, \textit{a posterior probability problem is completely determined if and only if it refers to a known phenomenon} (of which the characteristic function is known).

The problem of ``hypothesis likelihood''\footnoteM{See footnote \textit{a}.} consists of investigating the probability of a hypothesis or \textit{cause} to which a random phenomenon can be attributed, based on the observation of the frequency of successes of the phenomenon in a given number of trials. And we can conclude:\textit{ a problem of hypothesis likelihood is completely determined when and only when the phenomenon is known} (its characteristic function is given), \textit{the influence of the hypothesis is known} (the characteristic function of the subordinate phenomenon conditioned on the hypothesis is known), \textit{and the prior probability of the hypothesis itself is known}.

Otherwise, these two problems do not make sense.

The theorem in § 11 states everything which can be said precisely in a tentative inversion of Bernoulli's asymptotic theorem\footnoteM{Commonly referred to as the Law of Large Numbers.}.  As the number of trials grows, the probability of a random phenomenon tends to become equal to the frequency (with the restriction provided therein). But the convergence is not uniform for all the characteristic functions, and therefore, regardless of how large the number of trials carried out might be, it is not possible to deduce that the probability be approximately equal to the frequency without knowing in advance what the characteristic function of the phenomenon was. Nevertheless, we can say that as the number of trials grows, the conditions that have to be satisfied by the characteristic function of the phenomenon in order for the probability and frequency to be approximately equal become less restrictive.

These conclusions and examples can clarify the influence that empirical data\footnoteM{Lit.~\textit{dati dell'esperienza}, ``data of experience".} have on the evaluation of probabilities. They do so in a precise way for the two cases dealt with here, namely the posterior probability and the probability of causes, but also - in spirit - for the general case\footnoteF{A more exhaustive discussion is deferred to \textit{Probabilismo. Saggio critico sulla teoria delle probabilità e sul valore della scienza.} Biblioteca di Filosofia diretta da A.~ALIOTTA, Perrella ed., Napoli, 1931 (L.~5); see especially no.~22.}.

\section*{Acknowledgement}
We thank Valerio Varricchio and Fulvia de Finetti for their valuable comments regarding this translation. We are especially grateful to the latter for her permission to make this document publicly available for the benefit of the community.

\bibliography{references}
\bibliographystyle{plainnat}

\end{document}